\documentclass[12pt]{article}
\usepackage{amsfonts,amssymb,amsmath,epsfig,graphics,graphicx}

\usepackage[english]{babel}
\newcommand{\halmos}{\rule{1ex}{1.4ex}}

\makeatletter \@addtoreset{equation}{section} \makeatother

\newtheorem{ittheorem}{Theorem}
\newtheorem{itlemma}{Lemma}
\newtheorem{itproposition}{Proposition}
\newtheorem{itdefinition}{Definition}
\newtheorem{itcorollary}{Corollary}
\newtheorem{itremark}{Remark}
\newtheorem{itexamples}{Examples}

\newenvironment{theorem}{\addtocounter{equation}{1}
\begin{ittheorem}}{\end{ittheorem}}

\newenvironment{lemma}{\addtocounter{equation}{1}
\begin{itlemma}}{\end{itlemma}}

\newenvironment{proposition}{\addtocounter{equation}{1}
\begin{itproposition}}{\end{itproposition}}

\newenvironment{definition}{\addtocounter{equation}{1}
\begin{itdefinition}}{\end{itdefinition}}

\newenvironment{corollary}{\addtocounter{equation}{1}
\begin{itcorollary}}{\end{itcorollary}}

\newenvironment{remark}{\addtocounter{equation}{1}
\begin{itremark}}{\end{itremark}}

\newenvironment{examples}{\addtocounter{equation}{1}
\begin{itexamples}}{\end{itexamples}}

\newenvironment{proof}{\noindent {\em Proof}.\,\,\,}
{\hspace*{\fill}$\halmos$\medskip}

\newcommand{\beq}{\begin{eqnarray}}
\newcommand{\eeq}{\end{eqnarray}}

\newcommand{\beqq}{\begin{eqnarray*}}
\newcommand{\eeqq}{\end{eqnarray*}}

\newcommand{\bex}{\begin{examples}}
\newcommand{\eex}{\end{examples}}

\newcommand{\be}{\begin{equation}}
\newcommand{\ee}{\end{equation}}

\newcommand{\bl}{\begin{lemma}}
\newcommand{\el}{\end{lemma}}

\newcommand{\br}{\begin{remark}}
\newcommand{\er}{\end{remark}}

\newcommand{\bt}{\begin{theorem}}
\newcommand{\et}{\end{theorem}}

\newcommand{\bd}{\begin{definition}}
\newcommand{\ed}{\end{definition}}

\newcommand{\bp}{\begin{proposition}}
\newcommand{\ep}{\end{proposition}}

\newcommand{\bc}{\begin{corollary}}
\newcommand{\ec}{\end{corollary}}

\newcommand{\bpr}{\begin{proof}}
\newcommand{\epr}{\end{proof}}

\newcommand{\bi}{\begin{itemize}}
\newcommand{\ei}{\end{itemize}}

\newcommand{\ben}{\begin{enumerate}}
\newcommand{\een}{\end{enumerate}}

\newcommand{\Z}{\mathbb Z}
\newcommand{\R}{\mathbb R}
\newcommand{\N}{\mathbb N}

\newcommand{\pee}{\mathbb P}

\newcommand{\s}{\ensuremath{\mathcal{S}}}

\newcommand{\ii}{\ensuremath{\mathcal{I}}}

\newcommand{\al}{\ensuremath{\alpha}}

\newcommand{\ve}{\ensuremath{\varepsilon}}

\newcommand{\ga}{\ensuremath{\gamma}}
\newcommand{\Om}{\ensuremath{\Omega}}

\newcommand{\qaw}{N^{x,a}_t}
\newcommand{\qew}{N^{x,f}_t}
\newcommand{\da}{\dagger}

\begin{document}
\title{{\bf Freezing transitions in non-Fellerian particle systems}
}

\author{
C.\ Maes \footnote{Instituut voor Theoretische Fysica, K.U.Leuven,
Celestijnenlaan 200D, 3001 Leuven, Belgium.
email: Christian.Maes@fys.kuleuven.be}\\
F.\ Redig \footnote{Mathematisch Instituut, Universiteit Leiden,
Snellius, Niels Bohrweg 1, 2333 CA Leiden, The Netherlands.
\hskip 0.9cm
email: redig@math.leidenuniv.nl}\\
E.\ Saada \footnote{CNRS, UMR 6085, Laboratoire de
Math\'{e}matiques Rapha\"{e}l Salem, Universit\'{e} de Rouen,
Avenue de l'Universit\'{e}, BP.12, 76801 Saint-Etienne-du-Rouvray
Cedex, France.\hskip 3.6cm
email: Ellen.Saada@univ-rouen.fr}}

\maketitle

\footnotesize

\begin{quote}
{\bf Abstract}:  Non-Fellerian processes show phenomena that are
unseen in standard interacting particle systems.  We consider
freezing transitions in one-dimensional non-Fellerian processes
which are built from the abelian sandpile additions to which in
one case, spin flips are added, and in another case, the so
called anti-sandpile subtractions.  In the first case and as a function
of the sandpile addition rate, there is a sharp transition from a
non-trivial invariant measure to the invariant measure of the
sandpile process.  For the combination sandpile plus
anti-sandpile, there is a sharp transition from one frozen state
to the other anti-state
.
\end{quote}

\normalsize \vspace{12pt}
\section{Introduction}
Much of the motivation in the study of interacting particle
systems has come from the search for new phenomena. The abelian
sandpile model has been widely studied in the context of so called
self-organized criticality. From the point of view of probability
theory, it is the best known case of a spatially extended
non-Fellerian stochastic dynamics, \cite{M}.  It has challenged
our basic understanding of the construction of interacting processes in infinite volume,
even in one dimension.

In one dimension, the stationary measure of the standard abelian sandpile model is trivial in the
thermodynamic limit. However the dynamics of relaxation to this
measure is non-trivial. In \cite{MRSV} we have constructed the
dynamics for the one-dimensional sandpile model on the infinite
lattice $\Z$. The result is a monotone non-Fellerian process which
converges in finite time to its unique stationary state, which is
concentrating on the maximal configuration. As soon as one changes
to other lattices, such as decorated one-dimensional lattices, the
triviality of the limiting stationary measure disappears, and
especially in one dimension existence of thermodynamic limits is
not guaranteed due to the presence of infinite avalanches 
\cite{jarai2}.

In the present paper, we combine the one-dimensional sandpile
model with a spin-flip dynamics (pure spin flip as well as Glauber
type or more general spin flip processes with positive rates).
Indeed, in one dimension, the standard sandpile model has only two
possible heights per site, and spin flip just means changing the
height from one to the other. In the language of sandpiles, adding
a pure spin flip is the simplest example of combining two
different toppling mechanisms, the spin-flip part corresponding to
a purely dissipative (diagonal) toppling matrix. More precisely,
in Section 2, our dynamics has a formal generator
\[
Lf(\eta)= \alpha\sum_{x} [f(a_x \eta)-f(\eta)] + \sum_{x\in\Z}
c(x,\eta)\,[f(\theta_x\eta)-f(\eta)]
\]
where $a_x$ denote the abelian sandpile addition operators, and
$\theta_x$ the flip operator, on configurations $\eta \in
\{1,2\}^\Z$. In words that means that at rate $\alpha$ we add and
stabilize according to the abelian sandpile rule, and at rate
$c(x,\eta)$, we just flip the value of the height where we added.
In that way we have a parameter $\alpha$ that describes the
relative weight of the sand additions versus the spin flips. The
resulting ``sand-flip'' dynamics shows a freezing transition as a
function of that $\alpha$.   In the simplest case where
$c(x,\eta)=1$ (adding pure spin flip), our main result says that
for $\alpha\geq 1$ there is a finite time after which the system
reaches the maximal configuration (i.e., the sandpile part
``wins"), whereas the unique stationary measure is non-trivial and
mixing under spatial translations for $\alpha <1$. That is a
strong manifestation of the non-locality of the dynamics. Indeed,
for Fellerian processes such phenomenon cannot occur.

As a second example of such ``competition of different additions", in Section 3
we consider a combination of a sandpile and an
anti-sandpile process. This dynamics is inspired by \cite{km}. The
anti-sandpile part of the dynamics consists of removing grains and
stabilizing by reverse topplings. The infinite volume limit of the
anti-sandpile stationary measure is a Dirac measure concentrating
on the minimal configuration. Our main result here is that unless
the rates of addition and subtraction are equal, the limiting
stationary measure is a Dirac measure corresponding to the
dominant rate. We thus have a sharp transition between two
different frozen states.

These nonequilibrium phase transitions have an interest of their
own but they also go some way in adding extra and physically
relevant interactions to the standard abelian sandpile.  We have
in mind the sticky sandpiles of \cite{dharsticky}, for which our
dynamics is a subclass, and for which various transitions have been numerically checked.

As a final note, it is interesting to make the analogy with
non-Gibbsian measures.  In some sense, they are the ``equilibrium
analogue'' of non-Fellerian interacting particle systems.  In \cite{ms} an example of a freezing transition was
obtained, strongly connected with the absence of continuity of the
local conditional probabilities.

\section{Adding spin flips to the sandpile process}
The state space of our process is $\Om = \{1,2\}^{\Z}$. For a
configuration $\eta\in\Om$, $\eta(x)\in \{ 1,2\}$ is usually
interpreted as the height or the number of grains at site $x$.
That language can be continued even when combining the sandpile
automaton with other dynamics but we prefer to use the words
``active'' for $\eta(x)=2$ and ``inactive'' for $\eta(x)=1$.\\
The dynamics will change the configuration according to a
combination of the standard sandpile model and a spin flip
dynamics. We start with the simplest form of spin flip, 
changing ``active'' into ``inactive'' and vice versa at rate 1:
the spin flip $\theta_x$ is thus defined as 
\be\label{flip}
\theta_x\eta (y) =
\begin{cases}
(\eta(x) + 1)\mod 2,\ \text{if}\ y=x\\
\eta(y), \ \text{if}\ y\not= x
\end{cases}
\ee 
Only in Section \ref{gene} will we generalize the spin flip
part of the dynamics.\\
For the sandpile dynamics, we can rely on our previous work in
\cite{MRSV} where we have studied the infinite volume limit of the
one-dimensional sandpile process. We will therefore not bother to
redo the limiting procedures but below we immediately give the
result, the form of the infinite volume addition operators $a_x$.
The informal verbal prescription of the sandpile dynamics goes as
follows: if a site is inactive, it becomes active at rate
$\alpha$. If the site $x$ is already active, one looks left and
right of $x$ at the closest sites $x^-_\eta$ and $x^+_\eta$ which
are inactive. Again at rate $\alpha$ these two become active and
the mirror image of $x$ with respect to the middle of
$[x^-_\eta,x^+_\eta]$ becomes inactive.  That corresponds, in the infinite volume limit, to the
result (in finite volume) of adding and stabilizing through a sequence of topplings,
where upon a single toppling of a site the site looses two grains
and gives one grain to each neighbor, except if the site is at the
boundary where there is only one neighbor receiving a grain. See
\cite{MRSV} and \cite{redig} for more details on the abelian
sandpile model in $d=1$.

The infinite volume addition operator $a_x$ is defined more
precisely as follows:
 For $\eta\in\Om$ and
$x\in\Z$ with $\eta(x)=1$, we have $a_x\eta = \eta +e_x$ (where $e_x(x)=1$ and $e_x(y)=0$ otherwise),
i.e., inactive becomes active, or,  the height one at $x$ simply
changes to height two (and no other changes). For $\eta\in\Om$
and $x\in\Z$ with $\eta (x)=2$ we look at the right -- respectively
at the left -- of $x$ to find the first site $x^+_\eta$ (if that site
does not exist we put $x^+_\eta=\infty$), -- respectively $x^-_\eta$
(if that site does not exist we put $x^-_\eta=-\infty$) -- with $\eta(x^+_\eta)= \eta (x^-_\eta)=1$. We then define 
\be\label{ax}
a_x\eta (y) =
\begin{cases}
1 \ \text{if}\ y= x_\eta:=x^+_\eta+ x^-_\eta -x\\
2 \ \text{if}\ y\not=x_\eta,\ \text{and}\  x^-_{\eta}\leq y\leq x^+_\eta\\
\eta(y)\ \text{otherwise}
\end{cases}
\ee 
if both $x^+_\eta$ and $x^-_{\eta}$ exist.
In words, upon adding one unit at $x$, the first sites at
height one to the left ($x^-_\eta$) and to the right of $x$
($x^+_\eta$) become
 sites
with height $2$, and the site which is the mirror image of $x$
with respect to the middle of $x^+_\eta$, and $x^-_{\eta}$ becomes
of height one, all other sites remain unaltered. We have to extend
that definition to cases where one of the sites $x^+_\eta$,
$x^-_{\eta}$ does not exist (i.e., when there is no site to the
right or to the left of $x$ having height one). That is done by
taking the limit with ``boundary condition 1", i.e., if at least
one of the $x^{\pm}_\eta$ is infinite, then 
\be\label{inftycases}
a_x\eta (y) =
\begin{cases}
2 \ \text{if}\ 
x^-_{\eta}\leq y\leq x^+_\eta\\
\eta(y)\ \text{otherwise}
\end{cases}
\ee
Remark that \eqref{flip} is a special case of ``addition" with
``purely dissipative toppling", i.e., upon toppling an active 
site two grains disappear (diagonal toppling matrix).
 In that
sense combination of $a_x$ and $\theta_x$ is the simplest example
of combining two different toppling mechanisms (matrices) in one
process.

\subsection{Construction}\label{constr}

Intuitively, our process is governed by two independent collections of Poisson processes,
$N^{x,f}_t, N^{x,a}_t$,
indexed by sites $x\in\Z$, and independent for different sites. On the event
times of  $\qaw$ we apply the addition operator, and on the event times of
$\qew$ we ``flip" the state, i.e., we apply $\theta_x$. We put the rate of the
``sandpile-clocks" equal to $\alpha$, and the rate of the ``flip-clocks" equal to one.
Formally, our process has as a generator on local functions $f:\Om\to\R$,
\beq\label{geni}
Lf(\eta) &=& \alpha\sum_{x\in\Z} (f(a_x\eta)-f(\eta)) + \sum_{x\in\Z} (f(\theta_x\eta)-f(\eta))
\nonumber\\
&:=&
\alpha L_S f(\eta) + L_F (\eta)
\eeq
where $L_S$ stands for ``sandpile generator" and $L_F$ for ``flip-generator".

To show that there exists a Markov process with c\`adl\`ag-paths
corresponding to the Poisson process description above or to the
formal generator \eqref{geni}, we use a monotonicity argument
analogous to the one in \cite{MRSV}. We repeat the main steps and
the minor modifications to be done here. First we define the
action of $a_x$ on $\eta$ as a ``birth" if $\eta(x)=1$, and as an
avalanche if $\eta(x)=2$, whereas the (identical) action of
$\theta_x$ is (of course) also called a birth if $\eta(x)=1$, and a
death if $\eta(x)=2$. We can then split the formal generator in
three parts: \be\label{split} L= L_a + L_b + L_d \ee where
\beq\label{drie}
L_a f (\eta) &=& \alpha\sum_{x\in\Z} \chi(\eta(x) =2) (f(a_x\eta)-f(\eta))\nonumber\\
L_b f(\eta) &= & (1+\alpha)\sum_{x\in\Z} \chi(\eta(x) =1) (f(\theta_x\eta)-f(\eta))\nonumber\\
L_d f(\eta) &=& \sum_{x\in\Z} \chi(\eta(x)=1) (f(\theta_x\eta)-f(\eta))
\eeq
Here, ($\chi(.)$ denotes the indicator function).
The construction is then as follows:
\begin{itemize}
\item Construct a process corresponding to $L_a+L_d$ (only
avalanches and deaths) on the set $\Om_f$ of configurations
with a finite number of sites with height $2$. That is a
(non-explosive) countable state space Markov chain on the set of
finite subsets of $\Z$. Show by coupling that that process is {\em
monotone}. The coupling is identical to that of \cite{MRSV} for
the avalanche events. For the deaths: we let two ones die together
if possible, and otherwise independently. \item Construct a
process corresponding to $L_a+L_d$ with births in a finite
interval, i.e., having generator
\[
L_n f(\eta)= (L_a + L_d)f(\eta) + \sum_{x=-n}^n \chi(\eta(x)=1)
(f(\theta_x\eta)-f(\eta))
\]
We construct that process once more as a countable state space
Markov chain, and show that it is monotone. Its semigroup
$e^{tL_n} $ is denoted by $S_n(t)$. Moreover, we have the
following monotonicity as a function of the interval on which we
allow births: for all $t>0$, $n\in\N$, $f$ monotone,
$\eta\in\Om_f$,
\[
(S_n (t) f)(\eta) \leq (S_{n+1} (t) f)(\eta)
\]
\item
For general monotone $f$ and $\eta\in\Om$ arbitrary:
\be\label{semi}
S(t)f(\eta): = \sup_{n\in\N}\sup_{\eta\in\Om_f} S_n (t) f(\eta)
\ee
\end{itemize}
The process obtained by the above construction is called the
SF-process (sand-flip process). We denote its path space measure starting from $\eta$ by $P_{\eta}$.

\subsection{Basic properties}
Besides monotonicity, the SF-process has very similar
``quasi-Feller" properties as the one-dimensional sandpile process
of \cite{MRSV}. In particular, we have the following analogue of
Theorem 5.1 of \cite{MRSV}. Let us denote by
$\Om'\subset\Om$ the configurations with an infinite number
of ones to the left and to the right of the origin. We then
enumerate $\eta^{-1} \{1\} =\{ X_i(\eta), i\in\Z\}$ where
$X_0(\eta): = \min\{ x\geq 0:\eta(x)=1\}$, and the other $X_i$ are
in increasing order the sites where $\eta(x)=1$. The $X_i$ define
the $\eta$-dependent disjoint intervals $I_i = (X_{i-1} (\eta) ,
X_i (\eta)]$. A function is called $N$-local if it depends on the
heights $\eta(i)$ for $i\in \cup_{j=-N}^N I_j$. Every local
function is $N$-local, but a $N$-local function can be non-local,
e.g. $f(\eta)= e^{-|X_1 (\eta)|}$ is bounded $1$-local, but
non-local. The idea is that the natural space to define the action
of iterates of the generator is the set of $N$-local functions.
That is made precise in the following definition and theorem.

\bd\label{decco}
A configuration $\eta\in\Om'$ is called decent if
\be\label{dec} a(\eta) = \limsup_{n\to\infty} \frac{1}{2n}
\sum_{i=-n}^n |X_{i}(\eta)-X_{i-1}(\eta)| <\infty \ee The set of
decent configurations is denoted by $\Om_{dec}$. \ed

\begin{proposition}\label{seriesdef} Let $\eta\in\Om_{dec}$, $f$ be bounded
and $N$-local, then for $t<1/[4(1+\alpha)e a(\eta)]$, the
series $\sum_{n=0}^\infty [t^n(L^n f)(\eta)]/(n!)$ converges
absolutely and equals $S(t)f (\eta)$, where $S(t)$ is the
semigroup of the process defined above. In particular \be
\lim_{t\to 0} \frac{S(t) f(\eta)-f(\eta)}{t} = Lf(\eta) \ee
 i.e., $L$ is the ``pointwise generator" of the process.
 \end{proposition}

  \bpr The
same proof of \cite{MRSV} can be used, if one notices that the
extra ``death" part of the generator can only split one of the
intervals $I_i$ into smaller ones, by creating an extra $1$. This
implies that if $f$ is $N$-local, then $[f(a_i\eta)- f(\eta)]=0$
for all $i\in \Z\setminus \cup_{j=-N-1}^{N+1}I_j$. Therefore $Lf$
depends only on the heights in $\cup_{j=-N-1}^{N+1}I_j$. Iterating
the argument, one sees that $L^n f$ depends only on height in
$\cup_{j=-N-n}^{N+n}I_j$, and one recovers the same estimate
\be\label{crux} \|(L^n)f\|_\infty \leq  \prod_{k=0}^n \left(
\sum_{i=0}^{N+k} |I_i|\right) 2^n (1+\alpha)^n \|f\|_\infty \ee
which gives the result of the theorem, by application of lemma 4.1
in \cite{MRSV}. \epr

The following result, analogous to Corollary 6.1 in \cite{MRSV}
and to Proposition 3.1 in
\cite{MQ}, shows that the process is always non-Feller.
\begin{proposition}
 For all $\alpha >0$, the SF-process is non-Feller.
 \end{proposition}
 
 \bpr 
We denote by $\overline{2}$ the maximal configuration
$\eta\equiv 2$. We define the configuration $\eta_{spec}$  by 
\be \eta_{spec}(x) =
\begin{cases}
1 \ \text{if}\ x=0\\
2\ \text{otherwise}
\end{cases}
\ee
Then one shows as in \cite{MQ} that for $f_0 (\eta) =\eta(0)$
\be\label{imm}
\lim_{t\to 0, t>0} S(t) f_0 (\eta_{spec}) = 2
\ee
i.e., by the avalanche part of the dynamics, the isolated 1 is turned
``immediately" into a 2. Therefore,
the right limit of $\eta_t$ as $t\to 0, t>0$ is almost surely equal
to $\overline{2}$ when we start from $\eta_{spec}$.

This lack of right-continuity contradicts the Feller property. Indeed, if $S(t)$ were a
Feller-semigroup, then there would exist a uniformly dense set of
continuous functions which are in the domain of the generator,
i.e., for which 
\be\label{imma} 
\lim_{t\to 0, t>0} \frac{\|S(t)f -f\|_\infty}{t}=\|Lf\|_\infty 
\ee 
 However such a dense
set of continuous functions contains a function $f$ such that
\be\label{immo} 
f (\eta_{spec}) \not= f(\overline{2}) 
\ee
Combination of \eqref{imm}, \eqref{imma}, \eqref{immo} gives a
contradiction. 
\epr

\subsection{Stationary measure} 
Denote
by $\ii$ the set of invariant probability measures of the SF-process defined
in the previous section, by $\s$ the set of translation invariant probability measures on $\Om$. By monotonicity of the process $\ii$ is
non-empty. In fact we have \bt\label{SFthm} For all $\alpha >0$,
the SF-process is ergodic, i.e.,
$
\ii = \{ \mu_\alpha\}
$
and for all initial probability measures $\nu$ on $\Om$,
\be\label{ergosum} \lim_{t\to\infty} \nu S(t) = \mu_\alpha \ee
Moreover,
\begin{itemize}
\item for $\alpha <1$, the density of sites with height one is
given by \be \int \chi(\eta(0)=1)d\mu_\alpha =\frac{1-\alpha}{2} \ee
and $\mu_\alpha$ is a translation invariant measure which is
mixing under translations, non-product and gives positive measure
to all local events (i.e., has positive cylinders).

\item For $\alpha\geq 1$, \be \mu_\alpha = \delta_{2}
 \ee
  the
Dirac measure concentrating on the maximal configuration
$\eta\equiv 2$. Moreover for $t > [\log(\alpha+1) -
\log(\alpha-1)]/2$ and for every $\eta\in\Om$,
\[
\pee_\eta (\eta_t(0)=2)=1
\]
\end{itemize}
\et

\bpr We start with the following lemma. \bl\label{popo} Let $\mu$
be a probability measure on $\Om$ that is mixing under spatial
translations, with $\int \chi(\eta(0) =1) d\mu= \rho >0$. Then we have
\begin{itemize}
\item[a)] If $t<\rho/[4(1+\al)e]$, then $\mu S(t)$ is mixing under spatial translations. In particular,
\be\label{hypothese1.5}
\lim_{|x|\to\infty}
\int|S(t)[f\tau_x g] - S(t)f S(t)(\tau_x g)|d\mu=0
\ee
for all local functions $f,g$ on $\Om$, where $\tau_x$ denotes the spatial shift by $x\in\Z$.

\item[b)] Let $t'(\rho)$ be the solution of
\be\label{t'}
\Bigl(\rho+\frac{\al-1}{2}\Bigr)e^{-2t'(\rho)}+\frac{1-\al}{2}=0
\ee
when it exists, otherwise put
$t'(\rho)=+\infty$.
Then we have
for all $t< t'(\rho)$,
\beq\label{*}
\frac{d}{dt} \int\eta(0)d(\mu S(t))&=&\int L \eta (0)d(\mu S(t))\\
\label{***}
\rho(t):=\int \chi(\eta(0)=1)d(\mu S(t))&=&
\rho e^{-2t} + \frac{1-\alpha}{2} (1-e^{-2t})
\eeq
\end{itemize}
\el
\bpr
a) Is exactly \cite{MRSV}, Lemma 6.1.

b) Let $t_0=\inf\{\rho/5(1+\al)e, t'(\rho)\}$; by a), (\ref{*}) holds for
$t\leq t_0$. Denote \beq\label{KAAPLUS}
k^+ (i,\eta)&=& \inf \{j\geq 0:\eta (i+j) =1 \}\\
\label{KAAMIN}
k^- (i,\eta) &=& \inf \{ j>0: \eta (i-j) =1 \}
\eeq

We compute
\be\label{rene}
L\eta(0)=
\al\chi(\eta(0)=1)(k^+(1,\eta)+ k^-(0,\eta)+1)+3-\al-2\eta(0)
\ee
For $\nu\in\s$ (see \cite{MRSV}, (6.77)),
\[
\int  \chi(\eta(0) =1 )k^- (0,\eta )d\nu=
\int \chi(\eta(0)=1)(k^+ (1,\eta )+1)d\nu =1
\]
so that \be\label{int-L(eta(0))} \int L\eta(0) d(\mu S(t))=
\al+3-2\int \eta(0) d(\mu S(t)) =\al-1+2\int \chi(\eta(0)=1) d\mu
S(t) \ee and hence
\[
\frac{d\rho(t)}{dt} = -\alpha +1 - 2\rho (t)
\]
which gives \eqref{***} for $t<t_0$. If $t_0\not= t'(\rho)$, then
we can start the reasoning anew and iterate from $t_0$ to $t_1=
\min \{\rho(t_0)/[5(1+\al)e], t'(\rho (t_0))\}$, with new initial
distribution $\mu S(t_0)$. As a consequence, for $0\leq t\leq t_0
+t_1$, $\rho(t)$ is still given by \eqref{***}, etc. \epr

By monotonicity, the limits \be\label{exlim}
\nu_i=\lim_{t\to\infty}\delta_i S(t) \ee for $i=1,2$ exist and
define invariant measures with $\nu_1\leq \nu_2$. Moreover, the
process is totally ergodic (in the sense of \cite{L}, chapter 1, definition 1.9), i.e., $\ii$ is a singleton if and only
if $\nu_1=\nu_2$. In that case $\nu_1=\nu_2$ is also mixing under
spatial translations, see \cite{Andjel} Theorem 1.4 ii) (indeed,
the proof of the last theorem is a computation that does not
involve the Feller property of the considered Markov semi-group,
and relies on \eqref{hypothese1.5}).

Let $\lambda_\rho$ denote the translation invariant product
measure on $\Om$ with $\lambda_\rho(\eta(0)=1)=\rho$.
Then, using monotonicity and lemma \ref{popo}, for $t$ small
enough, \be\label{bhbh} \int \chi(\eta(0)=1)\ d(\delta_2 S(t)) \leq
\int \chi(\eta(0)=1)\ d\lambda_\rho S(t) =\rho e^{-2t} + \frac{1-\alpha}{2}
(1-e^{-2t}) \ee Since $\rho <1$ is arbitrary, we conclude
\be\label{klkl} \int \chi(\eta(0)=1)\ d(\delta_2 S(t)) \leq
\frac{1-\alpha}{2} (1-e^{-2t}) \ee

First consider $\alpha <1$. Starting from $\delta_1$, we have
$t'(\rho) =t'(0)=\infty$, and we obtain from lemma \ref{popo} \be
\lim_{t\to\infty} \int \chi(\eta(0)=1)\ d(\delta_1 S(t))
=\lim_{t\to\infty}\left(\rho e^{-2t} + \frac{1-\alpha}{2}
(1-e^{-2t})\right)=\frac{1-\alpha}{2} \ee Similarly, for $\rho$
close to one, $t'(\rho)= +\infty$, so that we can take the limit
$t\to\infty$ in \eqref{bhbh}, \eqref{klkl} to obtain, using
monotonicity once more, \be \frac{1-\alpha}{2}=\lim_{t\to\infty}
\int \chi(\eta (0)=1) d(\delta_1 S(t)) \leq \lim_{t\to\infty} \int
\chi(\eta (0)=1) d(\delta_2 S(t)) \leq \frac{1-\alpha}{2} \ee
Therefore, 
\be 
\nu_1 (\eta(0)=1)=\nu_2(\eta(0)=1)=\frac{1-\alpha}{2} 
\ee 
and combining that with
$\nu_1\leq \nu_2$, gives $\nu_1=\nu_2$.

To see that $\mu_\alpha$ is not a product measure, observe that
because $\mu_\alpha\in\s$, it can be product only if
$\mu_\alpha=\lambda_\rho$ with $\rho =(1-\alpha)/2$. Therefore, for all
$f$ local $\int Lf d\lambda_\rho=0$. Indeed $\lambda_\rho$ concentrates
in decent configurations, so we are allowed to use the generator
by proposition
\ref{seriesdef}.
Notice that in the sandpile part of the dynamics two or more neighboring
ones can never be created. Therefore, one  computes the action of
the sandpile part of the generator on the function $H_n(\eta) = \chi(\eta(1)=\ldots =\eta(n)=1)$ with $n\geq 2$,
which gives after integration over the product measure
\be
\int L_S H_nd\lambda_\rho=
- n\rho^n -2\sum_{i=0}^\infty i (1-\rho)^{i} \rho^{n+1}
=- n\rho^n -2\rho^{n-1}(1-\rho)
\ee
For the spin-flip part one has
\be
\int L_F (H_n) d\lambda_\rho
=-n\rho^n+n\rho^{n-1}(1-\rho)
\ee
Therefore, for the combined dynamics, the condition that
$\lambda_\rho$ is stationary leads to
\be\label{qazqazqaz}
\int (\alpha L_S + L_F) (H_n) d\lambda_\rho = 0=
\rho^{n-1} \left(-n\rho + n(1-\rho) -n\alpha\rho - 2\alpha (1-\rho)\right)
\ee
which gives (for $n\geq 2$)
\be\label{qazqaz}
\rho = \frac{n-2\alpha}{2n +(n-2)\alpha}
\ee
Notice that equation \eqref{qazqazqaz} is {\em not} valid for $n=1$, because a {\em single one}
{\em can} be created in the sandpile dynamics, so an extra term $(1-\rho)$ should
be added for $n=1$. Since \eqref{qazqaz} should be valid for all $n\geq 2$, we obtain a contradiction.
Hence, the invariant measure is indeed non-product.

For $\alpha=1$, \eqref{bhbh} gives
\[
\int \chi(\eta(0)=1) d(\delta_2 S(t)) =0
\]
Therefore $\delta_2 S(t)=\delta_2$. On the other hand,
\eqref{popo} gives
\[
\lim_{t\to\infty}\int \chi(\eta(0)=1) d(\delta_1 S(t)) =0
\]
Therefore $\nu_1=\delta_2$, hence we obtain
$\nu_1=\nu_2=\delta_2$.

Finally, consider $\alpha >1$. Then we have \be\label{wewe} t'(1)=
\frac12 \log\left(\frac{\alpha+1}{\alpha-1}\right) \ee Hence, for
all $\ve>0$ there exists $t_0 <t'(1)$ such that for all
$t>t_0$
\[
\int \chi(\eta(0)=1) d(\delta_1 S(t)) <\ve
\]
Therefore, $\nu_1=\lim_{t\to\infty} \delta_1 S(t) = \delta_2$, and
we conclude from the inequalities
$\nu_1=\delta_2 \geq \nu_2$, $\nu_2\leq \delta_2$
that $\nu_1=\nu_2=\delta_2$.

Moreover, we obtain that from any initial distribution $\nu$, the
limiting measure $\delta_2$ is reached in finite time $T_\nu \leq
t'(1)$. 
\epr 

\br For $\alpha\geq 1$, we have $\int Lg
d\mu_\alpha\not=0$ for non constant $g$, since for the sandpile
part $\int L_Sg d\mu_\alpha =0$, whereas for the flip part $\int
L_F g d\mu_\alpha\not=0$. Therefore, in that case the invariant
measure cannot be found by solving $\int Lf d\mu=0$ for $\mu$,
which gives another argument for the non-Fellerian character of
the SF-process. 
\er

\subsection{Generalization}\label{gene}
In this section we consider more general local perturbations of
the sandpile generator. We show that the freezing phenomenon,
i.e., having $\delta_2$ as unique invariant measure for $\alpha$
large enough, and a non-trivial invariant measure for $\alpha$
small persists.

More precisely, we consider a formal generator of the type
\be\label{bombom} L=\alpha L_S + L_G \ee
where $L_G$ is the
generator of a spinflip dynamics (i.e., with possibly
configuration dependent rates):

\be\label{bambam} L_Gf (\eta) = \sum_{x} c(x,\eta)
(f(\theta_x\eta)- f(\eta)) \ee where the flip-rates $c(x,\eta)$
are supposed to be translation invariant, local and bounded from
below. Therefore, \be\label{bornes} m\leq c(x,\eta)\leq M\ee for
some $0<m\leq M<\infty$ independent of $\eta$.

To define this process, we use a series expansion as in theorem
\ref{seriesdef} to {\em define} the semigroup. Remark that since
we do not assume that $L_G$ is the generator of a monotone
process, the semigroup cannot be constructed by monotonicity.
Instead, $S(t) f(\eta)$ is defined by the series expansion as long
as the configuration $\eta$ is decent, and contrary to the monotone case,
this cannot necessarily be extended to non-decent configurations
such as the maximal configuration $\overline{2}$.

We then have the following 
\bt 
Consider the process with formal
generator \eqref{bombom}. We have \ben \item For $\alpha <m$,
$\delta_2$ is not an invariant measure. In fact, if $\mu\in\s$ is
an invariant measure, then \be \mu(\eta(0)=1)\geq
\frac{m-\alpha}{2M}\ee \item If $\alpha
>M$,
then for all $\mu\in\s$, with $\mu (\eta(0)=1)>0$, $\mu
S(t)\to\delta_2$ as $t\to\infty$. Therefore, $\delta_2$ is the
only possible invariant measure. \een \et \bpr Let $\mu\in\s$ be
such that $\mu(\eta(0)=1)
>0$. We denote $\rho_t =\mu S(t)(\eta(0)=1)$. Then by the obvious
generalization of lemma \ref{popo}, we write \be
\frac{d\rho_t}{dt}= -\alpha -2\int \chi(\eta(0)=1) c(0,\eta) d(\mu S(t))
+ \int c(0,\eta) d(\mu S(t)) \ee Therefore \be\label{qweqwe} -\alpha
-2M\rho_t + m\leq \frac{d\rho_t}{dt}\leq -\alpha -2m\rho_t + M \ee
Hence, if $\alpha <m$ and $\rho_t < (m-\alpha)/2M$,
\[
\frac{d \rho_t}{dt} >0
\]
so there can be no invariant measure $\mu$ with $\rho=
\mu(\eta(0)=1) < (m-\alpha)/2M$. That proves the first item of the theorem. For
the second item, if $\alpha >M$ and if $\rho_t >0$, \eqref{qweqwe} gives
\[
\frac{d \rho_t}{dt} <0
\]
and hence there cannot be an invariant measure with $\rho=
\mu(\eta(0)=1) >0$. \epr

\br Even if $\mu S(t)$ converges to $\delta_2$ for all $\mu\in\s$
with $\mu(\eta(0)=1)
>0$, we cannot conclude that $\delta_2$ is an invariant measure,
because the process is not Feller. To see what can happen,
consider the following example. Starting from a configuration
$\eta\in\Om$, we flip a $1$ to $2$ at rate $1$, independently
for all lattice sites, and if the configuration is $\overline{2}$,
then we flip it at rate one to the minimal configuration $\eta\equiv 1$, denoted by $\overline{1}$. This is  a
non-Fellerian process with $\mu S(t)\to\delta_2$ for all
translation invariant measures $\mu\not=\delta_2$, but clearly,
$\delta_2$ is not invariant. In fact, the process has no invariant
measure. \er If $L_G$ is the generator of a monotone process, then
more precise results can be obtained. For that case we will stick
to an explicit example where once more an explicit closed equation
for the density can be obtained. More precisely, we consider the
flip rates \be\label{werwer} c(x,\eta) = 1-\ga f_x(\eta)
(f_{x-1} (\eta) +f_{x+1}(\eta))\ee where
\[
f_x (\eta) =1-2\chi(\eta(x)=1)
\]
These rates correspond to the standard Glauber choice for
\[
\ga = \frac12 \tanh (2\beta)\in \left[-\frac12,\frac12\right]
\]
where $\beta$ denotes the inverse temperature (here without any
meaning except for an effective coupling constant). We then have
the following analogue of theorem \ref{SFthm}. \bt For the process
with formal generator \eqref{bombom}, and rates \eqref{werwer} we
have \ben\item For $\alpha <\alpha_c= 1-2\ga$, there exists a
unique non-trivial invariant measure $\mu_\alpha$ with
\be\label{wawawa} \mu_\alpha (\eta(0)=1)=
\frac12\left(1-\frac{\alpha}{\alpha_c}\right) \ee \item For
$\alpha \geq \alpha_c$ , $\delta_2$ is the unique invariant
measure. \een \et \bpr Since the rates \eqref{werwer} satisfy
definition 2.1 of chapter 3 in \cite{L}, the process with
generator $L_G$ is monotone. Therefore, we can construct the
process with generator \eqref{bombom} by monotonicity as in
section \ref{constr}, where we replace the coupling for the birth
and death part by basic coupling. In particular the thus obtained
generalized SF-process is monotone. For $\mu\in\s$ with $\mu
(\eta(0)=1) >0$ we obtain \be \frac{d\rho_t}{dt} = -\alpha +
(1-2\rho_t) (1-2\ga) \ee and from that equation, combined with
monotonicity we can proceed as in the proof of theorem
\ref{SFthm}. \epr

\br
As one would expect intuitively, the critical
value $\alpha_c$ is decreasing in $\ga$, i.e., the freezing is
enhanced by stronger coupling.

Another simple choice in which we explicitly see the effect on
$\alpha_c$  is obtained by adding a bias to the spin
flip.  Then, \eqref{werwer} becomes
\[ c(x,\eta) =
1-\kappa f_x(\eta)= (1-\kappa)\chi(\eta(x) =2)+ (1+\kappa)\chi(\eta(x)=1)
\]
and a similar  calculation yields the same result but with a
critical value that now equals $\alpha_c = 1 - \kappa$.
\er

\section{ Adding ``anti-additions'' to the sandpile process} 

\subsection{ The anti-sandpile  model}
In words, the anti-sandpile process is a process where grains are removed from a configuration $\eta\in\Om$, and afterwards (if necessary) the configuration is stabilized instantaneously by {\em reversed topplings}. 

We first define the finite-volume process.
If after removing grains the height is zero at one or more sites $x\in [-N,N]$, then the configuration stabilizes by
a sequence of reversed topplings. Upon a reversed toppling of a site $x\in [-N,N]$
the site {\em gains} two grains and each of its neighbors (in $[-N,N]$) looses one grain.
This means that in a reversed toppling, the boundary sites act as {\em a source} (instead of a sink in
the ordinary toppling rule). The anti-addition operator $a^\dagger_x$ is then defined as the stable result
of the subtraction of one unit at site $x$ and performing reversed topplings until the configuration is {\em stable} (i.e. height everywhere $1$ or $2$) again.
\br
The anti-addition operator should not be confused with the inverse of the addition operator.
In fact, if $\eta$ is recurrent, {\em and} $a^\da_x\eta$ is recurrent, then
$a_x^\da\eta = a_x^{-1} (\eta)$, but $a_x^\da(\eta)$ need not be recurrent if $\eta $ is.
\er

In finite volume, the generator of the anti-sandpile process is given by
\be\label{antigen}
L^\da = \sum_{x=-N}^N (a_x^\da-I)
\ee
where $I$ denotes the identity operator.
Remark that $a_x^\da = \theta a_x\theta$ and
\be
L^\da = \Theta L \Theta
\ee
where $\Theta$ is ``global spinflip", i.e.,
\be
\Theta f(\eta) =f (\theta \eta)
\ee
with $\theta (\eta) (x) = (\eta (x)+1)\mod 2$.
Therefore the extension of the process generated by $L^\da$ to infinite volume is immediate. Its semigroup is given by
\be
S(t)^\da = \Theta S(t)\Theta
\ee
where $S(t)$ is the semigroup of the sandpile process.

In infinite volume, the
``anti-addition operator" $a^\da_x$ is then defined via
\[
a_x^\da \eta = \theta a_x\theta (\eta)
\]
Similarly, we introduce
\[
l^\pm (x,\eta)= k^\pm (x,\theta\eta)
\]
and the intervals $J_i(\eta) = I_i (\theta\eta)$.

\subsection{The SA process}

We now define the SA-process (i.e. ``sandpile + anti-sandpile'') as the process associated to the formal generator
\be\label{albet}
L_{\alpha\beta} =\alpha L + \beta L^\da
\ee
where $L= L_S = \sum_{x\in\Z} (a_x-I)$.

This process is constructed as follows: we {\em define} the semigroup acting on local functions
via the series expansion of proposition \ref{seriesdef}. 
This gives the finite dimensional distributions, and hence defines
a unique Markov process starting from decent configurations, where decent means here that {\em both} $\eta$ and $\theta\eta$ are decent in the sense of definition \ref{decco}.
We call the thus defined process the ``SA-process".

We then have the following.
\bp\label{SAmonotone}
The SA-process is monotone. As a consequence, it can be defined starting from any initial configuration.
\ep
\bpr
In \cite{MRSV} we constructed a generator for a coupling of the sandpile process which preserves the order.
The idea of this coupling is simply that if $\eta\leq\xi$, then for each site $j$ having height two in $\xi$ which
by an addition at some site $i\in\Z$ could be turned into a one, 
in $\eta$ either the height of $j$ is one or there exists a unique
site $x(j,\eta,\xi)$ having height two in $\eta$ such that addition at that site creates in $\eta$ a site of height one.

Let us call $L_S^c$ the (formal) generator of this coupling.
Remark now that if $\eta\leq\xi$ then of course $\Theta (\eta) \geq \Theta (\xi)$, and
for all $f$ monotone, $\Theta(f)$ is also monotone.

Therefore, the coupling with generator
\be\label{qqqq}
(L^c f)(\eta,\xi) = \alpha (L_S^c f)(\eta,\xi)  + \beta (L_S^c(\Theta f))(\theta\eta,\theta\xi)
\ee
defines a coupling that preserves the order.

This proves monotonicity. The consequence is clear since every configuration
can be written as an increasing (or decreasing) limit of decent configurations.
\epr

\subsection{Stationary measures for the SA process}

We then have the following analogue of theorem \ref{SFthm}.
\bt\label{SAthm}
Let $\ii$ be the set of invariant measures for the process with generator \eqref{albet}. Then we have
\begin{itemize}
\item
For $\alpha <\beta$
\[
\ii = \{\delta_1\}
\]
\item For $\alpha >\beta$,
\[
\ii =  \{ \delta_2\}
\]
\item For $\alpha=\beta$
\[
\ii_e\supset \{\delta_1,\delta_2\}
\]
\end{itemize}
\et
\bpr
We compute as before, starting from a translation invariant measure $\mu$ concentrating on decent configurations:
\[
\int L_S \chi(\eta(0)=1)d\mu  = -1 
\]
and hence since $\chi(\eta(0)=1)+ \chi(\eta(0)=2)=1$:
\[
\int L^\da \chi(\eta(0)=1) d\mu = +1
\]
Hence, starting from an initial measure $\mu$ on $\Om$ which is
translation invariant, mixing and concentrates on decent configurations, we obtain, using once more the notation 
$\rho_t = \int S(t)(\chi(\eta(0)=1))d\mu$:
\be\label{difeq} \frac{d\rho_t}{dt} = -\int (\alpha L + \beta
L^\da)(\chi(\eta(0)=1))d(\mu S(t))= (\beta-\alpha) \ee 
Of course this
equation is only valid as long as $0\leq (\beta-\alpha)t < 1$. It
expresses that the density of ones simply decreases or increases
linearly until no ones are present, resp.\ all sites are of height
one. 

Notice that we used here the analogue of proposition \ref{seriesdef}, which in this case implies that one
can use the generator as in the Feller case as long as acting on local functions and integrated over measures that have a non-zero density of sites having height two and of sites having height one.

Starting from this equation, one concludes that for $\alpha>\beta$
(and similarly for $\alpha <\beta$), there
can be no other invariant measure (which is also translation invariant) than $\delta_2$ (resp.\ $\delta_1$). 
We then deduce the first two statements of the theorem along the same lines as in theorem \ref{SFthm}, using monotonicity ( by proposition \ref{SAmonotone}).

For the last statement, use that if $\al=\beta$,
\[
\frac{d\rho_t}{dt} =0
\]
Therefore, using standard arguments based on monotonicity, 
we see that $\delta_1$ and $\delta_2$ are invariant measures.
\epr
\br
For the case $\al=\beta$ we have
\[
\frac{d\rho_t}{dt}=0
\]
i.e., the density is a conserved quantity. An open question here is whether in that
case for each density there exists a stationary
(in time) and ergodic (under translations) measure
with that density or whether the only extremal invariant measures
are $\{ \delta_1,\delta_2\}$. 

In that case, we can however say the following:
\ben
\item If $\mu$ is translation invariant, invariant for the dynamics, and with density $0<\rho<1$, then  $\mu\Theta$
is also invariant for the dynamics. Indeed, for any decent function $f$, $\Theta f$ is also decent, and $\int (L_S + \Theta L_S \Theta)f d\mu=0$ is equivalent to $\int (L_S + \Theta L_S \Theta)(\Theta f) d(\mu \Theta)=0$.
\item The product measure $\lambda_\rho$ with density 
$\lambda_\rho(\eta(0)=1)=\rho$ is not invariant. Indeed, we can proceed as in the proof of theorem \ref{SFthm}, and compute, for
$H_n(\eta) = \chi(\eta(1)=\ldots =\eta(n)=1)$ with $n\geq 2$,
\be\label{qaz}
\int (L_S + L^\dag) (H_n) d\lambda_\rho = \rho^{n-1} (1-\rho)(n - 2)
\ee
\een

\er

\noindent{\bf Acknowledgments.} E.S. thanks ESF, program RDSES,
EURANDOM and K.U.Leuven for financial support and hospitality.

\end{document}